\input amstex
\documentstyle{amsppt}
\magnification=\magstep1 \baselineskip=12pt \hsize=6truein
\vsize=8truein

\topmatter

\title
Higher order differential operators on manifolds without boundary
\endtitle
\author David Raske \\
\endauthor

\leftheadtext{Higher Order} \rightheadtext{David Raske}

\email david.t.raske{\@}gmail.edu \endemail

\abstract In this short note we review some results regarding the lack of maximum principles for higher order elliptic differential operators on manifolds without boundary. 
\endabstract

\endtopmatter

\document

\head 1. Higher order differential operators on manifolds without boundary \endhead

The aim of this note is to collect some results regarding the lack of maximum principles for higher order elliptic differential operators on manifolds. Most of these results are well known for elliptic operators defined on open subset of $\Bbb{R^n}$, with Dirichlet boundary conditions enforced, but there is little literature on what is the case for manifolds without boundaries. Since some positivity preserving takes place for higher-order differential operators with different boundary conditions enforced, one might naturally wonder if the absence of a boundary alters matters. As we will see this is not the case. 

For the sake of simplicity we will assume that our setting is a closed, smooth, Riemannian manifold. Let $P$ be a $2m$th, $m>1$ order elliptic differential operator. A classical result due to Calvert is the following

\proclaim{Theorem 1.1} Suppose $P$ is positive. Then the Green's function of $P+\lambda$ will not be positive for all $\lambda > 0$. \endproclaim

\demo{Proof} See [Ca1] and [Ca2]. \enddemo

For the case of elliptic operators on open subsets of $\Bbb{R^n}$ this is also proved in [CG]. Of course it follows that the hypothesis of positivity (and even self-adjointness) is not sufficient to guarantee the positivity of the Green's function in the manifold setting.

\proclaim{Theorem 1.2} Let $K$ be the heat kernel associated with $P$. Then $K$ must change signs for some $t>0$. \endproclaim

\demo{Proof} See [CG]. \enddemo

This result is somewhat surprising. One might expect the heat kernel associated with the bi-Laplacian to preserve positivity, but, as the above theorem notes, this is not the case. There is some hope, though, for pointwise lower bounds on heat kernels see [D]. 

As we will now see, the situation for eigenfunctions is similiary bleak.

\proclaim{Theorem 1.3} There exists a fourth order self-adjoint elliptic differential operator that does not have a signed principal eigenfunction nor is the principal eigenvalue simple. \endproclaim

\demo{Proof} See [HR]. \enddemo

For a interesting discussion of when positivty is guaranteed for the Poisson equation with positive data, or for the eigenproblem, the reader is advised to look at [R]. 

\vskip 0.1in \noindent {\bf References}:

\roster \vskip 0.1in \item"{[Ca1]}" B. Calvert, On T-accretive operators, Annali di Mat. Pura Appl. 94 (1972), 291-314.

\vskip 0.1in \item"{[Ca2]}" B. Calvert, Potential theoretic properties for accretive operators, Hiroshima Math J., 5 (No.3), (1975), 363-370.

\vskip 0.1in \item"{[CG]}" C, Coffman and C. Grover, Obtuse cones in Hilbert spaces and applications to partial differential equations, Jour. Funct. Anal., 35 (No.3), (1980), 369-396.

\vskip 0.1in \item"{[D]}" E.B.Davies, Pointwise lower bounds on the heat kernels of higher order elliptic operators, 125 (1999), 105-111.

\vskip 0.1in \item"{[HR]}"  E. Hebey and F. Robert,  Coercivity and Struwe's compactness for Paneitz operators with constant coefficients, Calc. Var., 13, (2001), 491-517.

\vskip0.1in \item"{[R]}" F. Robert, Fourth order equations with critical growth in Riemannian geometry, Preprint.   

\endroster

\enddocument